\documentclass[10pt]{amsart}
\usepackage{ amssymb,latexsym, amscd,pb-diagram}
\usepackage{sseq}
\usepackage[all]{xy}
\usepackage[square, numbers]{natbib}
\usepackage{graphicx}
\usepackage{mathrsfs}
\usepackage{color}
\usepackage{amsmath}
\usepackage{multirow}

\usepackage{tikz}
\usepackage{dsfont}
\usetikzlibrary{matrix}

 \newtheorem{theorem}{Theorem}[section]

 \theoremstyle{definition}
 
 \theoremstyle{definition}

 \numberwithin{equation}{section}

\newcommand{\ben}{\begin{equation}}
\newcommand{\een}{\end{equation}}

\newcommand{\integer}{\ensuremath{{\mathbb Z}}}

\newcommand{\complex}{\ensuremath{{\mathbb C}}}

\newcommand{\Cx}{\ensuremath{{\mathbb C}^*}}

\newcommand{\field}{\ensuremath{{\mathbb F}}}

\newcommand{\DD}{{\mathcal D}}

\newcommand{\VV}{{\mathcal V}}

\newcommand{\GG}{{\mathcal G}}
\newcommand{\CC}{\mathcal{C}}

\newcommand{\MM}{\mathcal{M}}
\newcommand{\HH}{\mathcal{H}}

\newcommand{\OO}{\mathcal{O}}

\newcommand{\Hom}{\mathrm{Hom}}

\newcommand{\Map}{\ensuremath{{\mathrm{Map}}}}

\newcommand{\Z}{{\mathbb Z}}

\newcounter{commentcounter}

\begin{document}

\title[Classification of Pointed Fusion Categories of dimension $p^3$]{Classification of Pointed Fusion Categories of dimension $p^3$ up to weak Morita Equivalence}

\thanks{The first author acknowledges the support 
of COLCIENCIAS through grant number FP44842-087-2017 of the Convocatoria Nacional J\'ovenes Investigadores e Innovadores No 761 de 2016. The third author acknowledges the financial support of the Max Planck Institute of Mathematics in Bonn, the Alexander von Humboldt Foundation, the Office of External Activities of the ICTP through Network NT8
 and COLCIENCIAS through grant number FP44842-013-2018 of the Fondo
Nacional de Financiamiento para la Ciencia, la Tecnolog\'ia y la Innovaci\'on.}
\author{Kevin Maya}
\address{Departamento de Matem\'{a}ticas y Estad\'istica, Universidad del Norte, Km.5 V\'ia Antigua a Puerto Colombia, 
Barranquilla, Colombia.}
\email{kjmaya@uninorte.edu.co}
\author{Adriana Mej\'ia Casta\~no}
\email{mejiala@uninorte.edu.co}
\author{Bernardo Uribe}
\email{bjongbloed@uninorte.edu.co}
\begin{abstract}
We give a complete classification of pointed fusion categories over $\complex$ of global dimension $p^3$ for $p$ any odd prime.
We proceed to classify the equivalence classes of pointed fusion categories of dimension $p^3$ and we  determine which of these equivalence classes have equivalent categories of modules. 
\end{abstract}

\maketitle

\section*{Introduction}
Fusion categories are rigid and  semisimple $\complex$-linear tensor categories whose isomorphism classes of simple objects are finitely many and the endomorphisms of the unit object is $\complex$. If all simple objects are invertible then it is called pointed fusion category. Pointed fusion categories are equivalent to fusion categories of the form $Vect(G, \omega)$ whose objects are 
complex vector spaces graded by the finite group $G$ and whose associativity constraint is encoded by a cocycle $\omega \in Z^3(G, \complex^*)$.

We say that two fusion categories $\DD$ and $\CC$ are weakly Morita equivalent if there exist
a module category $\MM$ of $\CC$ which is indecomposable and such that the tensor categories 
$\CC_\MM^*:=Fun_\CC(\MM,\MM)$ and $\DD$ are equivalent as tensor categories (see \cite[Def. 4.2]{MugerI}).

Naidu in \cite{Naidu} determined necessary and sufficient conditions  for two 
pointed fusion categories $Vect(G, \omega)$ and $Vect(\widehat{G}, \widehat{\omega})$ to be weakly Morita equivalent. 
The third author in \cite{Uribe} described these necessary and sufficient conditions with the use of the Lyndon-Hochschild-Serre spectral sequence
and gave an explicit description of the cocycles $\omega$ and $ \widehat{\omega}$ in order for the Morita equivalence to hold.

In this paper we follow the description done by the third author in \cite{Uribe}  in order to classify 
the Morita equivalence classes of pointed fusion categories of global dimension $p^3$ for $p$ odd prime. 
This description was used successfully in \cite{MunozUribe} where the classification
of pointed fusion categories up to Morita equivalence of global dimension $8$ was carried. This
paper is the continuation of \cite{MunozUribe} and finishes the classification for global dimension
$p^3$ for any prime.

This work is divided in two chapters. In the first chapter we setup explicit basis for $H^3(G, \complex^*)$ 
 and we calculate the space of orbits $H^3(G, \complex^*)/Aut(G)$ for each of the five groups of order $p^3$; this
 determines the equivalence classes of pointed fusion categories of global dimension $p^3$. In the second 
 chapter we recall the classification theorem of pointed fusion categories \cite[Thm. 3.9]{Uribe} and we 
 use the Lyndon-Hochschild-Serre spectral sequence associated to group extensions to explicitly find the Morita equivalence classes
 of pointed fusion categories of global dimension $p^3$. 
 
 Our results confirm the calculation of the number of pointed fusion categories of dimension $3^3$ up to Morita equivalence
 that was carried out in \cite[Page 34]{MignardSchauenburg} by computational methods, and permit us to determine the number of pointed fusion categories up to Morita equivalence
 of dimension $p^3$ for any fixed $p$.
 We obtain that there are
 $$5p+32$$
 Morita equivalence classes of pointed fusion categories of dimension $p^3$.
 
 Since this work is an application of the results obtained by the third author in \cite{Uribe} and \cite{MunozUribe}, we 
 will use the notation and the constructions done there. This work incorporates results that appear in the Master thesis 
 of the first author \cite{Maya}.

\section{Equivalence classes of pointed fusion categories of global dimension $p^3$}

 The pointed fusion categories $Vect(G, \omega)$ and $Vect(\widehat{G}, \widehat{\omega})$ are equivalent 
if and only if there is an isomorphism of groups $\phi: G \stackrel{\cong}{\to} \widehat{G}$ such that
$[\phi^* \widehat{\omega}]=[\omega]$ in $H^3(G, \complex^*)$. 

Therefore the equivalence classes of pointed fusion categories of global dimension $p^3$
can be determined as the union of the spaces of orbits $H^3(G,\complex^*)/Aut(G)$
where $G$ runs over the isomorphism classes of groups of order $p^3$.

The groups of order $p^3$ up to isomorphism are the abelian ones $\integer/p^3$,
$\integer/p^2 \times \integer/p$ and $(\integer/p)^3$ and the nonabelian ones
with presentations
$$\HH_p := \{ A,B,C | A^p=B^p=C^p=1, AC=CA, BC=CB, ABA^{-1}=BC\},$$
$$\GG_p:=\{a,b | a^p=b^{p^2}=1, aba^{-1}=b^{p+1}\}.$$
The first one $\HH_p$ is called the Heisenberg group of order $p^3$ and can be 
seen as the semi-direct product $(\integer/p \times \integer/p) \rtimes \integer/p$
with $\langle A \rangle $ acting on $\langle B,C \rangle$ leaving $C$ fixed and
sending $B \mapsto BC$. The second one $\GG_p$ 
can also be seen as the  semi-direct product $\integer/p^2 \rtimes \integer/p$
with $\langle a \rangle$ acting on $\langle b \rangle$ sending $b \mapsto b^{p+1}$.

We will make use of the ring structure of $H^*(G, \integer)$
in order to calculate $H^4(G, \integer) \cong H^3(G, \complex^*)$. In what
follows we will calculate $H^4(G, \integer)/Aut(G)$.

\subsection{The group $(\integer/p)^3$} \label{section (Zp)^3}
Let us recall some properties of the cohomology of $(\integer/p)^n$
that will be useful in what follows. The cohomology ring with 
coefficients in the field $\field_p$ of $p$-elements is
$$H^*((\integer/p)^n, \field_p) \cong \field_p[x_i,y_i | 1\leq i \leq n] /
 \langle x_i^2 | 1\leq i \leq n \rangle$$
 with $|x_i|=1$ and $|y_i|=2$. The Bockstein homomorphism
 $\beta: H^*((\integer/p)^n, \field_p) \to H^{*+1}((\integer/p)^n, \field_p)$
  is a graded derivation with $\beta(x_i)=y_i$ and $\beta(y_i)=0$.
  
  The integral cohomology ring is isomorphic to the kernel of the Bockstein
  $$H^*((\integer/p)^n, \integer) \cong \ker \left( \beta : H^*((\integer/p)^n, \field_p) \to H^{*+1}((\integer/p)^n, \field_p) \right)$$
and therefore we obtain
 $$H^4((\integer/p)^3, \integer) \cong\integer/p \langle y_1^2,y_2^2,y_3^2,y_1y_2, y_1y_3, y_2y_3, \beta(x_1x_2x_3) \rangle
\cong (\integer/p)^7.$$
The group of automorphisms of $(\integer/p)^3$ is isomorphic to $GL(3, \field_p)$ and for a matrix $A \in GL(3, \field_p)$
the action on the cohomology is given by its adjoint, i.e. $A^*x_i= A_{ij}x_j$ and $A^*y_i= A_{ij}y_j$. Note 
that $A^*\beta(x_1x_2x_3)= \det(A)\beta(x_1x_2x_3)$.

By the equivalence of quadratic forms and symmetric bilinear forms on a field $\field_p$ of odd prime characteristic,
we know that the action of the matrix $A$ on the quadratic forms $\integer/p \langle y_1^2,y_2^2,y_3^2,y_1y_2, y_1y_3, y_2y_3 \rangle$ is equivalent to the action $C \mapsto A^T C A$ on the set of symmetric $ 3\times 3$ matrices. In \cite[Thm. IV.10, p. 67]{Newman} it is shown that 
there are only two congruence classes of symmetric bilinear forms whenever the rank is fixed. Therefore the quadratic
forms on three variables are all congruent to one of the following 7 quadratic forms 
$$0, y_1^2, gy_1^2, y_1^2+y_2^2, gy_1^2+y_2^2, y_1^2+y_2^2+ y_3^2, gy_1^2+y_2^2+ y_3^2$$
where $g$ is any generator of the multiplicative group $(\integer/p)^*$.
Now, from the first 5 quadratic forms we may define the transformation that fixes $y_1$ and $y_2$ and maps
$y_3$ to $ay_3$ with $a \in (\integer/p)^*$. The determinant of this transformation is $a$ and therefore  we have the following $10$ orbits
in   $H^4((\integer/p)^3, \integer) / Aut((\integer/p)^3)$:
\begin{align*}
\{0\}, \OO(\beta(x_1x_2x_3)),  \OO(y_1^2), \OO(y_1^2 + \beta(x_1x_2x_3)), \OO(gy_1^2), \OO(gy_1^2 + \beta(x_1x_2x_3)),\\
\OO(y_1^2+y_2^2), \OO(y_1^2 +y_2^2 +\beta(x_1x_2x_3)), \OO(gy_1^2+y_2^2), \OO(gy_1^2 +y_2^2 +\beta(x_1x_2x_3)).
\end{align*}

If the matrix $A$ leaves the quadratic form $y_1^2+y_2^2+y_3^2$ fixed then $\det(A)^2=1$; we can choose $A=- \rm{Id}$ as one 
of those matrices. 
Therefore the remaining
$p+1$ orbits in $H^4((\integer/p)^3, \integer) / Aut((\integer/p)^3)$ are the following:
\begin{align*}
\OO(y_1^2+y_2^2+y_3^2 +a\beta(x_1x_2x_3)) \ \ \mbox{with} \ \ 0  \leq a \leq \frac{p-1}{2},  \\
\OO(gy_1^2+y_2^2+y_3^2 +a\beta(x_1x_2x_3)) \ \ \mbox{with} \ \ 0  \leq a \leq \frac{p-1}{2}.
\end{align*}

We conclude that the number of orbits is:
$$|H^4((\integer/p)^3, \integer) / Aut((\integer/p)^3)|=p+11.$$

\subsection{The group $\integer/p^2 \times \integer/p$} \label{subsection Zp2xZp} By K\"unneth's theorem we know that
$$H^4(\integer/p^2 \times \integer/p, \integer) \cong \integer \langle v^2, uv, u^2 \rangle /(p^2v^2,puv, pu^2) \cong
\integer/p^2 \oplus \integer/p \oplus \integer/p$$
where $H^*(\integer/p^2,\integer)=\integer[v]/(p^2v)$ and $H^*(\integer/p,\integer)=\integer[u]/(pu)$. 
Any automorphism of $\integer/p^2 \times \integer/p$ can be determined by the assignment
 \begin{align*}
\rho:\integer/p^2 \times \integer/p &\longrightarrow\integer/p^2 \times \integer/p\\
(1,0)&\longrightarrow (i,j)\\
(0,1)&\longrightarrow (pk,l)
\end{align*}
with $i \in (\integer/p^2)^*$, $j,k \in \integer/p$ and $l \in (\integer/p)^*$. The action of $\rho^*$ on $u$ and $v$
is $\rho^*u = lu+ pjv $ and $\rho^*v= ku+ iv$ and therefore
$$\rho^*u^2= l^2u^2, \ \ \rho^*uv= lku^2+iluv +pijv^2, \ \ \rho^*v^2=k^2u^2 + 2ikuv + i^2v^2.$$

By a direct calculation, we obtain 16 pairwise disjoint orbits
\begin{itemize}
    \item 4 orbits of size $\frac{p(p^2-p)(p-1)}{4}$: $$\OO(v^2+u^2),\OO(v^2+au^2), \OO(bv^2+u^2), \OO(bv^2+au^2), $$
    \item 4 orbits of size $\frac{p-1}{2}$: $$\OO(u^2),\OO(au^2),\OO(pv^2), \OO(bpv^2) $$
    \item 2 orbits of size $\frac{p(p^2-p)}{2}$: $$\OO(v^2),\OO(bv^2), $$
    \item 4 orbits of size $\frac{(p-1)^2}{4}$: $$\OO(pv^2+u^2),\OO(pv^2+au^2), \OO(bpv^2+u^2), \OO(bpv^2+au^2),$$
    \item 1 orbit of size $p^2(p-1)$: $$\OO(uv),$$
    \item 1 orbit of size $1$: $$\OO(0),$$
\end{itemize}
where $a\in (\integer/p)^*, b\in (\integer/p^2)^*$ such that $a$ and $b$ are not square numbers and $b\neq pk^2$ for any $k\in (\integer/p^2)^*$.
We conclude that the number of orbits is:
$$|H^4(\integer/p^2\times\integer/p, \integer) / Aut(\integer/p^2\times\integer/p)|=16.$$


\subsection{The cyclic group $\integer/p^3$} The multiplicative group $(\integer/p^3)^*$ of units of $\integer/p^3$
is isomorphic to the automorphism group $Aut(\integer/p^3)$  and acts on the generator of $H^*(\integer/p^3, \integer) = \integer[s]/(p^3s)$
by multiplication. Since the equation $i^2 \equiv j^2 \ \mbox{mod} \ p^n$ implies that $i\equiv j \ \mbox{mod} \ p^n$
or
$i \equiv -j \ \mbox{mod} \ p^n$, hence the action of $(\integer/p^3)^*$ on $H^4(\integer/p^3, \integer) = \langle s^2 \rangle/(p^3s)$
splits it into $7$ orbits: 
\begin{itemize}
    \item 2 orbits of equal size on the set $\{ks^2 | k \in (\integer/p^3)^*\}$: $$\OO(s^2),\OO(gs^2),$$
    \item 2 orbits of equal size on the set 
$\{jps^2 | j \in (\integer/p^2)^*\}$: $$\OO(ps^2),\OO(gps^2),$$
    \item 2 orbits of equal size on the set $\{lp^2s^2 | l \in (\integer/p)^*\}$: $$\OO(p^2s^2),\OO(gp^2s^2),$$
    \item the orbit of $0$,
\end{itemize}
where $g$ is any generator of the multiplicative group $(\integer/p)^*$.  We conclude that the number of orbits is:
$$| H^4(\integer/p^3, \integer)/ Aut(\integer/p^3)| = 7.$$


\subsection{The Heisenberg group $\HH_p$ } \label{subsection Heisenberg}

The automorphism group of the Heisenberg group 
$$\HH_p := \{ A,B,C | A^p=B^p=C^p=1, AC=CA, BC=CB, ABA^{-1}=BC\}$$ fits in the middle of the short exact
sequence
$$ 0 \to \integer/p \times \integer/p \to Aut(\HH_p) \to GL(2,p) \to 1 $$
where $\integer/p \times \integer/p$ corresponds to $Inn(\HH_p)$ and $GL(2,p)$ corresponds to $Out(\HH_p)$. 
The group $GL(2,p)$ is moreover isomorphic to the
automorphism group of  $\HH_p/Z(\HH_p)) \cong \integer/p \times \integer/p$ and the short exact
sequence splits since there a homomorphism
$$GL(2,p) \to Aut(\HH_p)$$ mapping the matrix $M=\left( \begin{matrix} a & b \\ c & d\end{matrix}\right)$
to the automorphism $\overline{M} \in Aut(\HH_p)$ with $\overline{M}(A)=A^aB^b$, $\overline{M}(B)=A^cB^d$
and $\overline{M}(C)=C^{ad-bd}=C^{det(M)}$. Hence the automorphism group is a semi-direct product
$$Aut(\HH_p) \cong (\integer/p \times \integer/p) \rtimes GL(2,p)$$
and the action of $Aut(\HH_p)$ on $H^*(\HH_p,\integer)$ can be calculated through the induced
action of $GL(2,p)$ on $H^*(\HH_p,\integer)$.

On order to determine the induced action of $GL(2,p)$ on $H^4(\HH_p,\integer)$ we will
use the Lyndon-Hochschild-Serre (LHS) spectral sequence associated to the short exact sequence
$$0 \to Z(\HH_p) \to \HH_p \to \integer/p \times \integer/p\to 0$$
where $Z(\HH_p) = \langle C \rangle \cong \integer/p$ is the center. 
Let the cohomology of the base be $H^*(\integer/p \times \integer/p, \field_p) \cong \field_p[w_1,w_2,z_1,z_2]/(w_1^2,w_2^2)$
with $\beta(w_i)=z_i$ and the cohomology of the fiber be $H^*(\integer/p,\integer)\cong \integer[t]/(pt)$.

The second page of the Lyndon-Hochschild-Serre spectral sequence 
$$E_2^{a,b}=H^a(\integer/p \times \integer/p, H^b(\integer/p, \integer))$$
has for relevant terms
\renewcommand{\arraystretch}{2}
$$\newcommand*{\tempb}{\multicolumn{1}{|c}{}}
\begin{array}{cccccccccc}
4 & \tempb & \langle t^2\rangle\\
3 & \tempb & 0& 0\\
2 & \tempb &\langle t\rangle&\langle tw_1, tw_2\rangle&\langle tz_1,tz_2,tw_1w_2\rangle\\
1 & \tempb & 0& 0 & 0 & 0\\
0 & \tempb & \integer&0& \langle z_1,z_2\rangle& \langle \beta(w_1,w_2)\rangle& \langle z_1^2,z_2^2,z_1z_2\rangle  \\ \cline{2-7}
 & & 0& ~~~~1~~~~ & ~~~~2~~~~ & ~~~~3~~~~ & 4
\end{array}
$$
and its isomorphic to the third page. The third differential $d_3$ maps $t$ to the Bockstein
of the $k$-invariant of the extension which in this case is the class $w_1w_2$, i.e.
$d_3(t)=\beta(w_1w_2)$. We have that $d_3(t^2)=2t\beta(w_1w_2)$, $d_3(tw_i)=\beta(w_1w_2w_i)=0$,
$d_3(tz_1)=\beta(w_1w_2)z_1$, $d_3(tz_2)=\beta(w_1w_2)z_2$ and $d_3(tw_1w_2)=0$. Hence
the fourth page of the spectral sequence has for relevant terms
\renewcommand{\arraystretch}{2}
$$\newcommand*{\tempb}{\multicolumn{1}{|c}{}}
\begin{array}{cccccccccc}
4 & \tempb & 0\\
3 & \tempb & 0& 0\\
2 & \tempb &0&\langle tw_1, tw_2\rangle&\langle tw_1w_2\rangle\\
1 & \tempb & 0& 0 & 0 & 0\\
0 & \tempb & \integer&0& \langle z_1,z_2\rangle& 0& \langle z_1^2,z_2^2,z_1z_2\rangle  \\ \cline{2-7}
 & & 0& ~~~~1~~~~ & ~~~~2~~~~ & ~~~~3~~~~ & 4
\end{array}
$$
thus obtaining that the associated graded $Gr(H^4(\HH_p,\integer))\cong \langle tw_1w_2,z_1^2,z_2^2,z_1z_2\rangle$. From \cite[Thm. 6.26]{Lewis} 
we know that $H^4(\HH_p,\integer)\cong \langle tw_1w_2,z_1^2,z_2^2,z_1z_2\rangle$
and therefore we may proceed to calculate the induced action of $GL(2,p)$ on
$H^4(\HH_p,\integer)$ with this choice of basis. For $M=\left( \begin{matrix} a & b \\ c & d\end{matrix}\right)$
the action of $\overline{M}$ is: $\overline{M}^*w_1=aw_1+cw_2$, $\overline{M}^*w_2=bw_1+dw_2$,
$\overline{M}^*z_1=az_1+cz_2$, $\overline{M}^*z_2=bz_1+dz_2$ and
$\overline{M}^*t=det(M)t$. Therefore 
\begin{align*}
\overline{M}^*(tw_1w_2)=det(M)^2tw_1w_2,  &&\overline{M}^*z^2_1=a^2z_1^2+2acz_1z_2 +c^2z_2^2,\\
\overline{M}^*z^2_2=b^2z_1^2+2bdz_1z_2 +d^2z_2^2, &&\overline{M}^*z_1z_2=abz_1^2+(ad+bc)z_1z_2 +cdz_2^2.
\end{align*}

For simplicity let us denote $\chi:=tw_1w_2$ and therefore
$$H^4(\HH_p,\integer) \cong \langle \chi, z_1^2,z_2^2,z_1z_2 \rangle \cong (\integer/p)^4.$$

By the classification of quadratic forms over a field $\field_p$ of prime odd characteristic
described in section \ref{section (Zp)^3} (cf. \cite[Thm. IV.10, p. 67]{Newman}), we know that there are $5$ orbits in
 $H^4(\HH_p,\integer)/Aut(\HH_p)$
of elements without the component $\chi$, these are:
$$\{0\}, \OO(z_1^2), \OO(gz_1^2), \OO(z_1^2+z_2^2), \OO(gz_1^2+z_2^2),$$
 where $g$ any generator of the multiplicative group $(\integer/p)^*$.
 
 If a matrix $A$ leaves $z_1^2+z_2^2$ or $gz_1^2+z_2^2$  fixed, then $\det(A)^2=1$ and therefore it acts trivially on $\chi$. 
 Hence we have $2(p-1)$ more orbits which are:
 \begin{align*}
 \OO(z_1^2+z_2^2 + a \chi) \ \ \mbox{with} \ \ 0 < a < p,\\
  \OO(gz_1^2+z_2^2 + a \chi) \ \ \mbox{with} \ \ 0 < a < p.
\end{align*}
The last $6$ orbits are
 \begin{align*}
 \OO(z_1^2 + \chi), \OO(z_1^2 + g\chi), \OO(gz_1^2 + \chi), \OO(gz_1^2 +g \chi), \OO(g\chi), \OO(\chi)
 \end{align*}
 where $g$ is any generator of the multiplicative group $(\integer/p)^*$. This follows from the fact that the matrix $A$
 that leaves $z_1^2$ fixed acts by multiplication by $\det(A)^2$ on $\chi$, and $\det(A)$ runs over all numbers in $(\integer/p)^*$.
 
We conclude that the number of orbits is: $$|H^4(\HH_p,\integer)/Aut(\HH_p)|= 2p+9.$$
 
\subsection{The group $\GG_p$}\label{section Gp} Any automorphism of the group
$$\GG_p:=\{a,b | a^p=b^{p^2}=1, aba^{-1}=b^{p+1}\}$$
is of the form
$$b \mapsto b^i a^j, \ \ a \mapsto b^{mp}a$$
with $i \in (\integer/p^2)^*$ and $j,m \in \integer/p$. Since the automorphisms $b \mapsto b$ and $a \mapsto b^{mp}a$
are all inner automorphism, we will only concentrate our attention on the ones of the form
$b \mapsto b^i a^j, \ \ a \mapsto a$. These automorphisms are generated by the automorphisms
$\rho(b)=b^i$, $\rho(a)=a$ and $\tau(b)=ba$, $\tau(a)=a$ with $i \in (\integer/p^2)^*$.

The LHS spectral sequence associated to the split extension
$$ 0 \to \langle b \rangle \to \GG_p \to \integer/p \to 0$$
has for second page $E_2^{n,m}\cong H^m(\integer/p , H^n(\integer/p^2,\integer))$
and its relevant terms are
\renewcommand{\arraystretch}{2}
$$\newcommand*{\tempb}{\multicolumn{1}{|c}{}}
\begin{array}{cccccccccc}
4 & \tempb & \langle pr^2 \rangle\\
3 & \tempb & 0& 0\\
2 & \tempb &\langle pr \rangle&0&0\\
1 & \tempb & 0& 0 & 0 & 0\\
0 & \tempb & \integer&0& \langle \gamma\rangle& 0& \langle \gamma^2 \rangle  \\ \cline{2-7}
 & & 0& ~~~~1~~~~ & ~~~~2~~~~ & ~~~~3~~~~ & 4
\end{array}
$$
where $H^*(\integer/p^2,\integer)=\integer[r]/(p^2r)$, $H^2(\integer/p^2,\integer)^{\integer/p}=\langle pr \rangle$,
$H^4(\integer/p^2,\integer)^{\integer/p}=\langle pr^2 \rangle$ and $H^*(\integer/p,\integer)=\integer[\gamma]/(p\gamma)$.

The automorphism $\rho$ induces an automorphism of $\integer/p$,
 and since the pullback $H^4(\HH_p,\integer) \to H^4(\integer/p^2,\integer)$ is injective 
 we conclude that 
 $$\rho^*(pr^2)= i^2pr^2, \ \ \rho^*(\gamma^2)=\gamma^2.$$
 
 The classes $pr$ and $\gamma$ represent explicit homomorphisms from $\GG_p$ to $\complex^*$;
 the first one sends $b \mapsto e^{2 \pi i /p}$ and $a \mapsto 1$ and the second sends $b \mapsto 1$ and
 $a \mapsto e^{2 \pi i /p}$. Therefore $\tau^* (pr)=pr$, $\tau^*\gamma = \gamma+pr$
 and we conclude that
 $$\tau^*(pr^2) = pr^2, \ \ \tau^*(\gamma^2)= \gamma^2.$$
 
Denoting the class $\delta:= pr^2$ we know that 
$$H^4(\GG_p,\integer) = \langle \delta, \gamma^2 \rangle \cong (\integer/p)^2$$
(cf. \cite[Thm. 5.2]{Lewis}), and therefore the action of $Aut(\GG_p)$ leaves the class $\gamma^2$ fixed,
and maps $\delta$ to $i^2\delta$ for $i \in (\integer/p^2)^*$. We conclude that for any prime
$p$ the orbits of 
$H^4(\GG_p,\integer)/Aut(\GG_p)$ are the following:
$$\{0\} , \{\gamma^2\}, \{2\gamma^2\}, ..., \{(p-1)\gamma^2\},$$
$$\OO(\delta), \OO(\gamma^2+\delta), \OO(2\gamma^2+\delta), ... , \OO((p-1)\gamma^2+\delta),$$
$$\OO(g\delta), \OO(\gamma^2+g\delta), \OO(2\gamma^2+g\delta), ... , \OO((p-1)\gamma^2+g\delta),$$
where $g$ is any generator of the multiplicative group $(\integer/p)^*$.

We conclude that the number of orbits is:
$$|H^4(\GG_p,\integer)/Aut(\GG_p)|=3p.$$

Collecting the last results, we obtain the next theorem:

\begin{theorem}\label{iso}
 There are $6p+43$ equivalence classes of pointed fusion categories of global dimension $p^3$. 
\end{theorem}\qed

\section{Morita equivalence classes of pointed fusion categories of global dimension $p^3$}

Now we will calculate which pointed fusion categories of global dimension $p^3$ have equivalent
categories of modules categories. We will base our calculations
on the classification theorem \cite[Thm. 3.9]{Uribe} and therefore we will use the same notation of \cite{Uribe}.
The following summary is taken from \cite[\S 2]{MunozUribe}

An skeletal indecomposable module category $\MM = (A \backslash G ,\mu)$ of
the skeletal category $\CC= \VV(G, \omega)$ of $Vect(G, \omega)$ is determined by a transitive $G$-set $K :=A \backslash  G$
with $ A $ subgroup of $G$, and a cochain 
$\mu \in C^2(G, \Map(K, \Cx))$ such that $\delta_H \mu = \pi^* \omega $ 
with $\pi^* \omega(k; h_1,h_2,h_3)= \omega(h_1,h_2,h_3)$ (see \cite[\S 3.3]{Uribe}).
The skeletal tensor category of the tensor category
 $\CC^*_\MM= Fun_\CC(\MM,\MM)$ is equivalent to one of 
the form $\VV(\widehat{G},\widehat{\omega})$ whenever $A$ is
 normal and abelian in $G$ \cite{Naidu} and if there exists
  a cochain $\gamma \in C^1(G, \Map(K, \Cx))$ such that 
  $\delta_G \gamma= \delta_K \mu$. 
  In particular this implies that the cohomology class of $\omega$
belongs to the subgroup of $H^3(G, \Cx)$ defined by
$$\Omega(G;A) := \ker \left( \ker \left( H^3(G, \Cx) \to E^{0,3}_\infty \right) \to E^{1,2}_\infty  \right),$$
which fits into the short exact sequence \cite[Cor. 3.2]{Uribe}
 $$0 \to E^{3,0}_\infty \to \Omega(G;A) \to E^{2,1}_\infty \to 0$$
where $E_n^{*,*}$ denotes the $n$-th page of the Lyndon-Hochschild-Serre spectral sequence associated to the group extension
$1 \to  A \to G \to K \to 1$.

 Denote the dual group ${{\mathbb{A}}} := \Hom(A, \complex^*)$ and consider 
cocycles $F \in Z^2(K, A)$ and $\widehat{F} \in Z^2(K, {{\mathbb{A}}})$. Denote by $H= A \rtimes_F K$ and
$\widehat{H}= K \ltimes_{\widehat{F}} {{\mathbb{A}}}$ the groups defined by the multiplication laws 
$$(a_1,k_1) (a_2,k_2) := (a_1 ({}^{k_1}a_2) F(k_1,k_2),k_1k_2),$$
$$(k_1, \rho_1) \cdot (k_2, \rho_2) := (k_1k_2, (\rho_1^{k_2}) \rho_2 \widehat{F}(k_1,k_2))$$
respectively.  The necessary and sufficient  conditions for two pointed fusion categories to be Morita equivalent are the following (cf. \cite{Naidu}):

\begin{theorem}\cite[Thm. 5.9]{Uribe}  \label{classification theorem}
 Let $G$ and $\widehat{G}$ be finite groups, $\omega \in Z^3(G, \complex^*)$ and $\widehat{\omega} \in Z^3(\widehat{G}, \complex^*)$. Then the tensor categories $Vect(G,\omega)$ and $Vect(\widehat{G}, \widehat{\omega})$ are weakly Morita equivalent if and only if the following
  conditions are satisfied:
 \begin{itemize}
 \item There exist isomorphisms of groups 
 $$\phi : H= A \rtimes_F K \stackrel{\cong}{\to} G \ \ \ \ 
 \widehat{\phi} : \widehat{H}= K \ltimes_{\widehat{F}} {{\mathbb{A}}} \stackrel{\cong}{\to} \widehat{G}$$
 for some finite group $K$ acting on the abelian group $A$,
 with cocycles $F \in Z^2(K, A)$ and $\widehat{F} \in Z^2(K, {{\mathbb{A}}})$. 
 \item There exists $\epsilon : K^3 \to \Cx$ such that $\widehat{F} \wedge F =\delta_K \epsilon$.
 \item The cohomology classes satisfy the equations $[\eta]=[ \phi^* \omega ]$ and 
 $[\widehat{\eta}]=[\widehat{\phi}^*\widehat{\omega}]$ with 
 \begin{align*}
   \eta((a_1,k_1),(a_2,k_2),(a_3,k_3)) := & \widehat{F}(k_1,k_2)(a_3) \ \epsilon(k_1,k_2,k_3),\\
    \widehat{\eta}((k_1, \rho_1),  (k_2,\rho_2),(k_3 ,\rho_3)) :=& \epsilon(k_1,k_2,k_3) \ \rho_1(F(k_2,k_3)).
   \end{align*}
\end{itemize}
\end{theorem}\qed

The abelian groups ${{\mathbb{A}}}$ and $A$ are (non-canonically) isomorphic as $K$-modules, and hence
both $G$ and $\widehat{G}$ could be seen as extensions of $K$ by $A$. In order to calculate
all possible Morita equivalences, we will analyze the Morita equivalences that appear while fixing the 
group $K$ and the $K$-module $A$.

Let us recall the equivalence classes of normal abelian subgroups of the groups of order $p^3$. Two subgroups will be equivalent if there is an automorphism of the group that maps one to the other.
The following table contains the information the equivalence classes of these subgroups:
\begin{center}
 \begin{tabular}{|c|c|c|c|c|c|}
\hline
Isomorphic to&$\integer/p^3$ & $\integer/{p^2}\oplus\integer/{p}$ & $\HH_p$ & $\GG_p$\tabularnewline
\hline
 $\integer/p$ &$\left\langle p^2 \right\rangle $ & $\left\langle (p,0)\right\rangle $ & $\left\langle C\right\rangle $ & $\left\langle b^p\right\rangle $\tabularnewline
\hline
 $\integer/p$ & & $\left\langle (mp,1)\right\rangle, m \in (\integer/p)^* $
 & &
\tabularnewline
\hline
$\integer/p\oplus \integer/p$ &  &$\left\langle (p,0)\right\rangle \oplus\left\langle (0,1)\right\rangle $ & $\left\langle B,C \right\rangle $ & $\left\langle a,b^p \right\rangle $
\tabularnewline
\hline
$\integer/p^2$&$\left\langle p\right\rangle $ & $\left\langle (1,k)\right\rangle,  k \in \integer/p $ &  &$\left\langle ba^l\right\rangle, l \in \integer/p $
\tabularnewline
\hline
\end{tabular}
\end{center}
In  $(\integer/p)^3$ all subgroups of order $p$ are equivalent to $\langle (0,0,1) \rangle$ and all subgroups of order $p^2$ are equivalent to $\langle (0,1,0),  (0,0,1) \rangle$.

The procedure that we will follow is the same one that appeared in \cite[\S 2]{MunozUribe}, let us recall it.
 We fix the groups $K$ and $A$, we  take the groups
that are extensions of $K$ by $A$ and we take explicit choices of subgroups
from the table above that provide the extensions. Then we  calculate the relevant terms of the second page of the Lyndon-Hochschild-Serre
spectral sequence, which are the same for all extensions of $K$ by $A$ previously chosen, and we calculate the third page for each extension $0 \to A \to G \to K\to 1$. Then we determine the cohomology class
of the 2-cocycle $F$ that makes $G \cong A \rtimes_F K$ and we calculate the cohomology classes in $\Omega(G;A)$.
With this information and Theorem \ref{classification theorem} we determine the Morita equivalence classes of pointed fusion categories for groups 
that are extensions of $K$ by $A$.

\subsection{$K= \integer/p$ and $A=\integer/p^2$ with trivial action}

The two possible extensions are $\integer/8$ and $\integer/4 \times \integer/2$ with the following choices of subgroups:
\begin{align*} 
1 \longrightarrow\langle (1,0)\rangle \longrightarrow \integer/p^2 & \times \integer/p \longrightarrow \integer/p \longrightarrow 1  \\
1 \longrightarrow\langle p\rangle \longrightarrow & \integer/p^3  \longrightarrow \integer/p  \longrightarrow 1.
\end{align*}
For the group $\integer/p^2 \times \integer/p$ the relevant terms of the second page of the LHS spectral sequence are:
\renewcommand{\arraystretch}{2}
\begin{align*}\label{rob}
\newcommand*{\tempb}{\multicolumn{1}{|c}{}}
\begin{array}{cccccccccc}
4 & \tempb & \integer/p^2=\langle v^2\rangle\\
3 & \tempb & 0& 0\\
2 & \tempb & \Z/p^2& \Z/p & \mathbb{Z}/p =\langle uv\rangle\\
1 & \tempb & 0& 0&0 &0\\
0 & \tempb & \integer& 0 & \mathbb{Z}/p & 0 & \mathbb{Z}/p =\langle u^2\rangle\\ \cline{2-7}
& & 0~~~& 1~~~ & 2~~~ & 3~~~ & 4 
\end{array}
\end{align*}
Since $H^2(\integer/p^2, \integer/p) = \integer /p$ we have that 
$\integer/p^3 \cong \integer/p \ltimes_{uv} {\integer/p^2}$, and therefore the relevant terms of the fourth page of the LHS spectral sequence
for $\integer/p^3$ are:

\renewcommand{\arraystretch}{2}
\begin{align*}
\newcommand*{\tempb}{\multicolumn{1}{|c}{}}
\begin{array}{cccccccccc}
4 & \tempb & \integer/p^2=\langle s^2\rangle/(p^2s^2)\\
3 & \tempb & 0& 0\\
2 & \tempb & \Z/p^2& 0 & \mathbb{Z}/p =\langle p^2s^2\rangle\\
1 & \tempb & 0& 0&0 &0\\
0 & \tempb & \integer& 0 & \mathbb{Z}/p & 0 & 0\\
 \cline{2-7}
& & 0~~~& 1~~~ & 2~~~ & 3~~~ & 4 
\end{array}
\end{align*}
where in this case $d_3:E_2^{1,2} \stackrel{\cong}{\to} E_2^{4,0}$.

We get the split short exact sequences 
\begin{align*}0 \to  \langle u^2\rangle  &\to \Omega(\integer/p^2 \times \integer/p; \langle (1,0)\rangle) \to \langle uv\rangle \to 0\\
& 0 \to \Omega(\integer/p^3; \langle p^2 \rangle) \stackrel{\cong}{\to} \langle p^2s^2 \rangle \to 0
 \end{align*} and we conclude that the only Morita equivalence that appear 
is:
\begin{align*}
Vect(\integer/p^3,0)\simeq_M Vect(\integer/p^2\times\integer/p,uv).
\end{align*}
Note that the calculations of section \ref{subsection Zp2xZp} show that the classes
$kuv+lu^2$ 
with $k \in (\integer/p)^*$ and $l \in \integer/p$ belong to the same orbit
under the action of the automorphism group.

\subsection{$K=\integer/p$ and $A=\integer/p^2$ with non-trivial action.}
In this case the only possible extension is the group $\GG_p$, and since the groups 
$E_2^{2,2}=0 =E_2^{3,1}$ for the LHS spectral sequence, we conclude that there are
no non-trivial Morita equivalences associated to this case.

\subsection{$K=\integer/p^2$ and $A=\integer/p$.}
The two extensions are:
\begin{align*}
1  \longrightarrow\langle (0,1)& \rangle \longrightarrow \integer/p^2 \times \integer/p \to \integer/p^2 \longrightarrow 1\\
1  \longrightarrow&\langle p^2 \rangle \longrightarrow \integer/p^3 \longrightarrow \integer/p^2  \longrightarrow 1,
\end{align*}
and the relevant terms of the second page of the LHS spectral sequence for $\integer/p^2 \times \integer/p$ are:
\renewcommand{\arraystretch}{2}
$$\newcommand*{\tempb}{\multicolumn{1}{|c}{}}
\begin{array}{cccccccccc}
4 & \tempb & \mathbb{Z}/p=\langle u^2\rangle\\
3 & \tempb & 0& 0\\
2 & \tempb &\mathbb{Z}/p& \mathbb{Z}/p& \mathbb{Z}/p=\langle uv\rangle\\
1 & \tempb & 0& 0&0 &0\\
0 & \tempb & \integer& 0& \mathbb{Z}/p^2 & 0 & \mathbb{Z}/p^2=\langle v^2\rangle  \\ \cline{2-7}
 & & 0& ~~~~1~~~~ & ~~~~2~~~~ & ~~~~3~~~~ & 4
\end{array}
$$
We have that  $\integer/p^2=\integer/p^2 \rtimes_{uv}\integer/p$, and the relevant terms of the fourth page of the LHS spectral sequence
associated to the extension of $\integer/p^3$ are:
\renewcommand{\arraystretch}{2}
$$\newcommand*{\tempb}{\multicolumn{1}{|c}{}}
\begin{array}{cccccccccc}
4 & \tempb & \mathbb{Z}/p=\langle s^2\rangle/(ps^2)\\\
3 & \tempb & 0& 0\\
2 & \tempb &\mathbb{Z}/p& 0 & \mathbb{Z}/p=\langle ps^2\rangle/(p^2s^2)\\
1 & \tempb & 0& 0&0 &0\\
0 & \tempb & \mathbb{C}^*& \mathbb{Z}/p^2 & 0 & \mathbb{Z}/p=\langle p^2s^2\rangle  \\ \cline{2-7}
 & & 0& ~~~~1~~~~ & ~~~~2~~~~ & ~~~~3~~~~ & 4
\end{array}
$$
where in this case $d_3:E_3^{1,2}\cong \integer/p \to E_3^{4,0} \cong \integer/p^2$ is injective and therefore
$E_4^{4,0} \cong \integer/p$.

We obtain the short exact sequences 
\begin{align*}
&0 \to  \langle v^2 \rangle  \to \Omega(\integer/p^2 \times \integer/p; \langle (0,1) \rangle) \to \langle uv \rangle \to 0\\
 0 &\to \langle p^2s^2 \rangle /(p^3s^2)\to \Omega(\integer/p^3; \langle p \rangle) \to \langle ps^2 \rangle /(p^2s^2) \to 0
 \end{align*} 
 where $\Omega(\integer/p^3; \langle p \rangle) \cong \integer/p^2$. We conclude that the only Morita equivalences that appear 
are:

\begin{align*}
Vect(\integer/p^3,0)\simeq_M & Vect(\integer/p^2\times \integer/p,uv),\\
Vect(\integer/p^3,kp^2s^2)\simeq_M & Vect(\integer/p^2\times\integer/p,uv+kv^2) 
\end{align*} 
for $k \in (\integer/p^2)^*$.

By the calculations of section \ref{subsection Zp2xZp} we know that in the orbit
of $uv+kv^2$ appears $luv + i^2kv^2$ for any $l \in (\integer/p)^*$ and any $i \in (\integer/p^2)^*$.
Therefore the orbit of the class $p^2s^2$ in $H^4(\integer/p^3,\integer)$ is related to the orbit of the class $uv+v^2$ 
in $H^4(\integer/p^2 \times \integer/p, \integer)$ via a Morita equivalence. The same applies to the
complementary orbit of $p^2s^2$ on the set of classes $\{ms^2 | m \in (\integer/p)^* \}$.

\subsection{$K=\integer/p$ and $A= \integer/p \times \integer/p$ with trivial action}
The relevant extensions are:
\begin{align*}
1  \longrightarrow\langle (0,0,1), (0,1,0) \rangle \longrightarrow (\integer/p)^3 \to \integer/p  \longrightarrow 1\\
1  \longrightarrow\langle (p,0),(0,1) \rangle \longrightarrow \integer/p^2 \times \integer/p \longrightarrow \integer/p  \longrightarrow 1.
\end{align*}
The relevant terms for the second page $E_2^{a,b}=H^a(\integer/p,H^b(\integer/p \times \integer/p, \integer))$ of the LHS spectral sequence are:
\renewcommand{\arraystretch}{2}
\begin{align*}
\newcommand*{\tempb}{\multicolumn{1}{|c}{}}
\begin{array}{ccccccccccc}
4 & \tempb & \langle y_2^2,y_2y_3,y_3^2 \rangle\\
3 & \tempb & \langle \beta(x_2x_3) \rangle&\langle x_1\beta(x_2x_3) \rangle\\
2 & \tempb & \langle y_2,y_3 \rangle &\langle x_1y_2,x_1y_3 \rangle &\langle y_1y_2,y_1y_3 \rangle& \langle x_1y_1y_2,x_1y_1y_3 \rangle\\
1  & \tempb & 0 & 0 & 0 & 0 &0\\
0 & \tempb & \integer &0 & \langle y_1 \rangle & 0 & \langle y_1^2 \rangle& 0 \\ \cline{2-8}
& & 0& 1& 2 & 3 & 4 & 5
\end{array}
\end{align*}

The second differential of the LHS spectral sequence associated to the group $\integer/p^2 \times \integer/p$ is generated by
$$d_2x_2=y_1, \ d_2y_2=0, \ d_2x_3=0, \ d_2y_3=0$$
where the class $y_1$ is the $k$-invariant of the extension.
We calculate
$$d_2(\beta(x_2x_3))=-y_1y_3, \ d_2(x_1\beta(x_2x_3))=x_1y_1y_3  $$

and  the relevant terms of the third page of the spectral sequence become
\renewcommand{\arraystretch}{2}
\begin{align*}
\newcommand*{\tempb}{\multicolumn{1}{|c}{}}
\begin{array}{ccccccccccc}
4 & \tempb & \langle y_2^2,y_2y_3,y_3^2 \rangle\\
3 & \tempb & 0 &0 \\
2 & \tempb & \langle y_2,y_3 \rangle &\langle x_1y_2,x_1y_3 \rangle &\langle y_1y_2 \rangle& \langle x_1y_1y_2 \rangle\\
1  & \tempb & 0 & 0 & 0 & 0 &0\\
0 & \tempb & \integer &0 & \langle y_1 \rangle & 0 &  \langle y_1^2 \rangle  & 0 \\ \cline{2-8}
& & 0& 1& 2 & 3 & 4 & 5
\end{array}
\end{align*}

The only non-trivial third differential is
$$d_3(x_1y_2)= y_1^2$$
and the relevant terms of the fourth page of the spectral sequence become 
\renewcommand{\arraystretch}{2}
\begin{align*}
\newcommand*{\tempb}{\multicolumn{1}{|c}{}}
\begin{array}{ccccccccccc}
4 & \tempb & \langle y_2^2,y_2y_3,y_3^2 \rangle\\
3 & \tempb & 0 &0 \\
2 & \tempb & \langle y_2,y_3 \rangle &\langle x_1y_3 \rangle &\langle y_1y_2 \rangle\\
1  & \tempb & 0 & 0 & 0 & 0\\
0 & \tempb & \integer &0 & \langle y_1 \rangle & 0 & 0  \\ \cline{2-7}
& & 0& 1& 2 & 3 & 4 
\end{array}
\end{align*}

The correspondence of the previous terms with respect to the ones defined in \S \ref{subsection Zp2xZp} is the following: $\langle y_1y_2 \rangle$ corresponds
to the group $\langle p v^2 \rangle$, $\langle y_2^2 \rangle$ 
corresponds
to the group $\langle  v^2 \rangle/(pv^2)$, $\langle y_2y_3 \rangle$ corresponds
to the group $\langle  uv \rangle$ and $\langle y_3^2 \rangle$ corresponds
to the group $\langle u^2 \rangle$.

We get the split short exact sequences 
\begin{align*}0 \to  &\langle y_1^2 \rangle  \to \Omega((\integer/2)^3; \langle (0,0,1),(0,1,0) \rangle) \to \langle y_1y_2, y_1y_3 \rangle \to 0\\
 &0 \to  \Omega(\integer/p^2 \times \integer/p; \langle (p,0),(0,1) \rangle) \stackrel{\cong}{\to} \langle pv^2  \rangle \to 0
 \end{align*} 
 where all the non-trivial classes in $\langle y_1y_2, y_1y_3 \rangle$ define an extension isomorphic to $\integer/p^2 \times \integer/p$.
We conclude that the only Morita equivalence that appear is:
\begin{align*}
Vect(\integer/p^2 \times\integer/p ,0)\simeq_M Vect((\integer/p)^3,y_1y_2).
\end{align*}

\subsection{$K=\integer/p\times \integer/p$ and $A=\integer/p$}
In this case the relevant extensions are:
\begin{align*}
1 \longrightarrow \langle (0,0,1) \rangle  \longrightarrow& (\integer/p)^3 \longrightarrow \integer/p \times \integer/p \longrightarrow 1\\
1 \longrightarrow \langle (p,0) \rangle \longrightarrow \integer/p^2 &\times \integer/p \longrightarrow  \integer/p \times \integer/p \longrightarrow 1\\
1 \longrightarrow \langle C \rangle \longrightarrow & \HH_p \longrightarrow\integer/p \times \integer/p \longrightarrow 1\\
1 \longrightarrow \langle b^p \rangle \longrightarrow &\GG_p \longrightarrow \integer/p \times \integer/p \longrightarrow 1.
\end{align*}

The relevant terms of the second page $E_2^{a,b}=H^a(\integer/p \times \integer/p, H^b(\integer/p, \integer))$
of the LHS spectral sequence are:
\renewcommand{\arraystretch}{2}
$$\newcommand*{\tempb}{\multicolumn{1}{|c}{}}
\begin{array}{cccccccccc}
4 & \tempb & \langle y_3^2\rangle\\
3 & \tempb & 0& 0 & 0 \\
2 & \tempb &\langle y_3 \rangle& \langle y_3x_1, y_3x_2\rangle&\langle y_3y_1, y_3 y_2, y_3x_1x_2\rangle\\
1 & \tempb & 0& 0 & 0 &0 &0 \\
0 & \tempb &\integer & 0 & \langle z_1, z_2\rangle & \langle \beta(x_1x_2) \rangle&\langle z_1^2,z_1z_2,z_2^2 \rangle  \\ \cline{2-7}
 & & 0& ~~~~1~~~~ & ~~~~2~~~~ & ~~~~3~~~~ & 4
\end{array}
$$

The third differential of the spectral sequence incorporates the information of the $k$-invariant
of the extension $$\kappa \in H^2(\integer/p \times \integer/p, \integer/p) \cong \integer/p \langle y_1,y_2,x_1x_2 \rangle.$$
For the group $\integer/p^2 \times \integer/p$ we take $\kappa = y_1$, for $\HH_p$ we take $\kappa=x_1x_2$ and for
$\GG_p$ we take $\kappa=y_2+x_1x_2$. The third differential from the second row to the 0th-row
$$d_3 : H^a((\integer/p)^2; \integer/p) \otimes H^2(\integer/p; \integer/p) \to H^{a+3}((\integer/p)^2, \integer)$$
is given by the formula
$$d_3(y_3 p(x_1,x_2,y_1,y_2)) = \beta(\kappa \ p(x_1,x_2,y_1y_2))$$
where $p(x_1,x_2,y_1,y_2)$ is any polynomial. The third differential on $y_3^2$ is simply $d_3(y^2_3)= 2y_3\beta(\kappa)$.

\subsubsection{$(\integer/p)^3$} For the group $(\integer/p)^3$ we obtain a split extension
$$0 \to \langle y_1^2, y_2, y_1y_2\rangle \to \Omega((\integer/p)^3; \langle (0,0,1) \rangle) \to \langle y_3y_1,y_3y_2,y_3x_1x_2\rangle \to 0$$
where the classes $ly_1 +my_2$ with $lm \neq 0$ define an extension isomorphic to $\integer/p^2 \times \integer/p$,
the classes $kx_1x_2$ with $k \neq 0$ define an extension isomorphic to $\HH_p$ and the classes $kx_1x_2 +ly_1+my_2$ with $k\neq0$ and
$lm \neq 0$ define an extension isomorphic to $\GG_p$.
Therefore it is enough to analyze the cases determined by the classes $y_1$, $x_1x_2$ and $y_2+x_1x_2$. 

\subsubsection{$\integer/p^2 \times \integer/p$} For the group $\integer/p^2 \times \integer/p$ we have that it is isomorphic to the group
 $\integer/p\rtimes_{y_1}(\integer/p \times\integer/p)$ and the relevant terms of the fourth page of the LHS are:

\renewcommand{\arraystretch}{2}
\begin{align*}
\newcommand*{\tempb}{\multicolumn{1}{|c}{}}
\begin{array}{cccccccccc}
4 & \tempb & \langle y_3^2\rangle\\
3 & \tempb & 0& 0\\
2 & \tempb & \langle y_3 \rangle & 0& \langle y_3y_1,y_3y_2\rangle \\
1 & \tempb & 0& 0 & 0\\
0 & \tempb & \integer & 0&  \langle y_1,y_2\rangle & \langle \beta(x_1x_2) \rangle & \langle y_2^2 \rangle\\ \cline{2-7}
& & 0~~~& 1~~~ & 2~~~ & 3~~~ & 4
\end{array}
\end{align*}
Following the notations of section \ref{subsection Zp2xZp} we have that $\langle y_2^2 \rangle$ corresponds
to $\langle u^2 \rangle$,  $\langle y_3y_1 \rangle$ corresponds
to $\langle uv \rangle$, $\langle y_3y_2 \rangle$ corresponds
to $\langle pv^2 \rangle$ and $\langle y_3^2 \rangle$ corresponds
to $\langle v^2 \rangle/(pv^2)$. Therefore we obtain
$$0 \to \langle u^2 \rangle \to \Omega(\integer/p^2\times \integer/p; \langle (p,0) \rangle) \to \langle uv, pv^2 \rangle \to 0.$$
Since all the classes in $E_4^{2,2}$ induce extensions
isomorphic to $\integer/p^2 \times \integer/p$ we conclude
that the only Morita equivalences that we obtain in this case are:
\begin{align*}
Vect(\integer/p^2 \times \integer/p, 0) &\simeq_M Vect((\integer/2)^3,y_3y_1),\\
Vect(\integer/p^2 \times \integer/p, ku^2) & \simeq_M Vect((\integer/2)^3,y_3y_1 + ky_2^2)
\end{align*}
where $k \in (\integer/p)^*$. By the calculations of section \ref{subsection Zp2xZp} we know that the orbit
of $u^2$ in $H^4(\integer/p^2 \times \integer/p, \integer)$ is related via a Morita equivalence to the orbit of the class 
$y_3y_1 + y_2^2$ in $H^4((\integer/p^3, \integer)$ and if the complementary orbit of $u^2$ in the set $\{ k u^2 | k \in (\integer/p)^* \}$
is generated by $ju^2$, then its orbit is related to the orbit of the class $y_3y_1 + jy_2^2$.

\subsubsection{$\HH_p$} The calculation of the relevant terms of the third page of the LHS spectral sequence for the group 
$\HH_p \cong \integer/p \rtimes_{x_1x_2}(\integer/p \times \integer/p )$ was done in section \ref{subsection Heisenberg}
and therefore we obtain the short exact sequence
$$0 \to \langle z_1^2,z_1z_2,z_2^2 \rangle \to \Omega(\HH_p; \langle C \rangle) \to \langle \chi \rangle \to 0.$$
Since $E_4^{2,2}$ is generated by the class $y_3x_1x_2$ and $x_1x_2$ is the $k$-invariant for $\HH_p$, we only get
the following Morita equivalences:
\begin{align*}
Vect(\HH_p, &0) \simeq_M Vect((\integer/2)^3,\beta(x_1x_2x_3)),\\
Vect(\HH_p, az_1^2+bz_1z_2+cz_2^2&) \simeq_M Vect((\integer/2)^3,\beta(x_1x_2x_3)+ ay_1^2+by_1y_2+cy_2^2)
\end{align*}
for any $a,b,c \in \integer/p$.

\subsubsection{$\GG_p$} The relevant terms of the fourth page for the LHS spectral sequence for the group
$\GG_p \cong \integer/p \rtimes_{x_1x_2+y_2}(\integer/p \times \integer/p )$ are:
\renewcommand{\arraystretch}{2}
\begin{align*}
\newcommand*{\tempb}{\multicolumn{1}{|c}{}}
\begin{array}{cccccccccc}
4 & \tempb & 0\\
3 & \tempb & 0& 0\\
2 & \tempb & 0 & 0 &\langle y_3(y_2-x_1x_2)\rangle \\
1 & \tempb & 0& 0& 0 & 0\\
0 & \tempb & \integer & 0&\langle y_1,y_2\rangle & 0 &\langle y_1^2 \rangle \\ \cline{2-7}
& & 0~~~& 1~~~ & 2~~~ & 3~~~ & 4
\end{array}
\end{align*}
Following the notation of section \ref{section Gp} we have that $\langle y_3(y_2-x_1x_2)\rangle$ corresponds to $\langle \delta \rangle$
and $\langle y_1^2 \rangle$ corresponds to $\langle \gamma^2 \rangle$.
Therefore we obtain the split short exact sequence
$$0 \to \langle \gamma^2 \rangle \to \Omega(\GG_p; \langle b^p \rangle) \to \langle \delta \rangle \to 0.$$
Since the the class $y_2-x_1x_2$ is also a $k$-invariant of $\GG_p$, we only get 
the following Morita equivalences:
\begin{align*}
Vect(\GG_p, 0)& \simeq_M Vect((\integer/2)^3, y_3y_2-\beta(x_1x_2x_3)),\\
Vect(\GG_p, k\gamma^2) \simeq_M & Vect((\integer/p)^3,y_3y_2-\beta(x_1x_2x_3)+ ky_1^2)\\
\end{align*}
for $k \in (\integer/p)^*$. Note that the classes $y_3y_2-\beta(x_1x_2x_3)+ ky_1^2$ and $y_3y_2-\beta(x_1x_2x_3)+ ly_1^2$
for $k \neq l$ are in different orbits in $H^4((\integer/p)^3,\integer)$, therefore compatible with the fact that
the automorphism group of $\GG_p$ acts trivially on $\langle \gamma^2 \rangle$.

\subsection{$K=\integer/p$ and $A=\integer/p \times \integer/p$ with non-trivial action.}

Let us consider the action of $\integer/p$ on $\integer/p \times \integer/p$ generated by the assignment
$$(1,0) \mapsto (1,1), \ \ (0,1) \mapsto (0,1).$$
This action is compatible with the action of $A$ on $B$ since $ABA^{-1}=BC$ and with the action
of $b^{-1}$ on $a$ since $b^{-1}ab=ab^p$.
The non-isomorphic extensions are
\begin{align*}
1  \longrightarrow \langle B,C\rangle\longrightarrow \HH_p\longrightarrow \integer/p  \longrightarrow 1\\
1  \longrightarrow \langle a, b^p\rangle\longrightarrow \GG_p\longrightarrow \integer/p  \longrightarrow 1
\end{align*}
where $\HH_p \cong \langle B,C\rangle \rtimes \integer/p$ and $\GG_p$ has for $k$-invariant
any of the non-trivial elements in $H^2(\integer/p, \integer/p \times \integer/p) \cong \integer/p$.

  For $\HH_p$ the relevant elements of the second page of the LHS spectral sequence are
  $E_2^{r,s}\cong H^r(\integer/p, H^s(\integer/p \times \integer/p, \integer))$:
  
  \renewcommand{\arraystretch}{2}
\begin{align*}
\newcommand*{\tempb}{\multicolumn{1}{|c}{}}
\begin{array}{cccccccccc}
4 & \tempb & \integer/p \\
3 & \tempb & \integer/p& \integer/p\\
2 & \tempb & \integer/p & \integer/p &\integer/p \\
1 & \tempb & 0& 0& 0 & 0\\
0 & \tempb & \integer & 0& \integer/p & 0 &\integer/p \\ \cline{2-7}
& & 0~~~& 1~~~ & 2~~~ & 3~~~ & 4
\end{array}
\end{align*}
  where $E_2^{4,0}= \integer/p$ corresponds to $\langle z_1\rangle $ and $E_2^{2,2}= \integer/p$ corresponds to $\langle z_1z_2\rangle$.
  
  For $\GG_p$ the relevant terms of the third page of the LHS spectral sequence are:
  \renewcommand{\arraystretch}{2}
\begin{align*}
\newcommand*{\tempb}{\multicolumn{1}{|c}{}}
\begin{array}{cccccccccc}
4 & \tempb & \integer/p \\
3 & \tempb & 0& \integer/p\\
2 & \tempb & \integer/p & \integer/p &0 \\
1 & \tempb & 0& 0& 0 & 0\\
0 & \tempb & \integer & 0& \integer/p & 0 &\integer/p \\ \cline{2-7}
& & 0~~~& 1~~~ & 2~~~ & 3~~~ & 4
\end{array}
\end{align*}
and the relevant terms of the fourth page are:
    \renewcommand{\arraystretch}{2}
\begin{align*}
\newcommand*{\tempb}{\multicolumn{1}{|c}{}}
\begin{array}{cccccccccc}
4 & \tempb & \integer/p \\
3 & \tempb & 0& \integer/p\\
2 & \tempb & \integer/p &0 &0 \\
1 & \tempb & 0& 0& 0 & 0\\
0 & \tempb & \integer & 0& \integer/p & 0 &0 \\ \cline{2-7}
& & 0~~~& 1~~~ & 2~~~ & 3~~~ & 4
\end{array}
\end{align*}
where $E_4^{0,4}$ corresponds to $\langle \gamma^2 \rangle$ and $E_4^{1,3}$ corresponds to $\langle \delta\rangle$.

Therefore we get
\begin{align*}0 \to \langle z_1^2  \rangle  \to \Omega(\HH_p; \langle B, C\rangle) \to \langle z_1z_2 \rangle \to 0 \ \ \mbox{and} \ \  \Omega(\GG_p; \langle a, b^p\rangle) = 0
 \end{align*} 
and we conclude that the only Morita equivalence that appear are:
   \begin{align*}
   Vect(\GG_p, 0) \simeq_M Vect(\HH_p,z_1z_2+lz_1^2)  \end{align*} 
for any  $l \in \integer/p$.

\subsection{} We conclude that the only weak Morita equivalences between pointed fusion categories
of global dimension $p^3$ are the ones that appear on each row of the following table:
\begin{center}
\scalebox{0.7}{
\begin{tabular}{ ||c|c|c|c|c|| } 
 \hline 
 $\integer/p^3$ & $\integer/p^2 \times \integer/p$ & $(\integer/p)^3$ &  $\HH_p$ & $\GG_p$\\[0.5ex] 
 \hline\hline
 $\{0\}$ & $\OO(uv)$ &&& \\ \hline
  $\OO(p^2s^2)$ & $\OO(uv+v^2)$ &&& \\ \hline
   $\OO(gp^2s^2)$ & $\OO(uv+gv^2)$ &&& \\ \hline
   & $\{0\}$ & $\OO(y_1y_2)$&& \\ \hline
    & $\OO(u^2)$ & $\OO(y_1y_2+y_2^2)$&& \\ \hline
     & $\OO(gu^2)$ & $\OO(y_1y_2+gy_2^2)$&& \\ \hline
     && $\OO(\beta(x_1x_2x_3))$& $\{0\}$& \\ \hline
  && $\OO(\beta(x_1x_2x_3) + y_1^2)$& $\OO(z_1^2)$& \\ \hline
  && $\OO(\beta(x_1x_2x_3) + gy_1^2)$& $\OO(gz_1^2)$& \\ \hline
  && $\OO(\beta(x_1x_2x_3) + hy_1^2 + y_2^2)$& $\OO(hz_1^2 + z_2^2)$& \\ \hline
    && $\OO(y_2y_3+ \beta(x_1x_2x_3))$& $\OO(z_1z_2)$&$\{0\}$ \\ \hline
   && $\OO(y_2y_3+ \beta(x_1x_2x_3) + y_1^2)$& &$\{\gamma^2 \}$ \\ \hline
    && $\OO(y_2y_3+ \beta(x_1x_2x_3) + 2y_1^2)$& &$\{2\gamma^2 \}$ \\ \hline
     && \vdots & &\vdots \\ \hline
      && $\OO(y_2y_3+ \beta(x_1x_2x_3) + (p-1)y_1^2)$& &$\{(p-1)\gamma^2 \}$ \\ \hline
\end{tabular}}
\end{center}
where $g$ is any generator of the multiplicative group $(\integer/p)^*$. Here the quadratic form $hz_1^2 + z_2^2$ is either
$z_1^2 + z_2^2$ or $gy_1^2 + y_2^2$ depending on which one is not congruent to $z_1z_2$; for example
for $p=3$ the form $z_1+z_2$ is not congruent to $z_1z_2$, and for $p=5$ the form
$2z_1+z_2$ is not congruent to $z_1z_2$.

Hence we conclude that there are $p+10$
 non-trivial Morita equivalence classes of pointed fusion categories of dimension $p^3$, $p+9$ with
 two classes and only one
 with three classes.
 
Since there are $6p+43$ equivalence classes of pointed fusion categories of dimension $p^3$ by Theorem \ref{iso}, subtracting $p+11$ classes that are Morita equivalent to others, we obtain the following result:
 \begin{theorem}
 There are $5p+32$
 Morita equivalence classes of pointed fusion categories of dimension $p^3$.
 \end{theorem}
 
 For $p=3$ therefore there are $61$
    equivalence classes of pointed fusion categories of dimension $27$ and $47$ Morita equivalence classes
    of the same dimension.
    The following table provides an explicit description in terms of cohomology classes of the list that appears in \cite[Page 34]{MignardSchauenburg}
and confirms the calculations done in GAP for $p=3$:
 
\begin{center}
\scalebox{0.7}{
\begin{tabular}{ ||c|c|c|c|c|| } 
 \hline 
 $\integer/27$ & $\integer/9 \times \integer/3$ & $(\integer/3)^3$ &  $\HH_3$ & $\GG_3$\\[0.5ex] 
 \hline\hline
 $\{0\}$ & $\OO(uv)$ &&& \\ \hline
  $\{9s^2\}$ & $\OO(uv+v^2)$ &&& \\ \hline
   $\{18s^2\}$ & $\OO(uv+2v^2)$ &&& \\ \hline
   & $\{0\}$ & $\OO(y_1y_2)$&& \\ \hline
    & $\OO(u^2)$ & $\OO(y_1y_2+y_2^2)$&& \\ \hline
     & $\OO(2u^2)$ & $\OO(y_1y_2+2y_2^2)$&& \\ \hline
     && $\OO(\beta(x_1x_2x_3))$& $\{0\}$& \\ \hline
  && $\OO(\beta(x_1x_2x_3) + y_1^2)$& $\OO(z_1^2)$& \\ \hline
   && $\OO(\beta(x_1x_2x_3) + 2y_1^2)$& $\OO(2z_1^2)$& \\ \hline
     && $\OO(\beta(x_1x_2x_3) + y_1^2 +y_2^2)$& $\OO(z_1^2+z_2^2)$& \\ \hline
    && $\OO(y_2y_3+ \beta(x_1x_2x_3))$& $\OO(z_1z_2)$&$\{0\}$ \\ \hline
   && $\OO(y_2y_3+ \beta(x_1x_2x_3) + y_1^2)$& &$\{\gamma^2 \}$ \\ \hline
    && $\OO(y_2y_3+ \beta(x_1x_2x_3) + 2y_1^2)$& &$\{2\gamma^2 \}$ \\ \hline
\end{tabular}}
\end{center}

For $p=5$ there are $73$
    equivalence classes of pointed fusion categories of dimension $125$ and $57$ Morita equivalence classes of the same
  dimension.
    The following table gives the explicit Morita equivalence classes in terms of cohomology classes:
\begin{center}
\scalebox{0.7}{
\begin{tabular}{ ||c|c|c|c|c|| } 
 \hline 
 $\integer/27$ & $\integer/9 \times \integer/3$ & $(\integer/3)^3$ &  $\HH_3$ & $\GG_3$\\[0.5ex] 
 \hline\hline
 $\{0\}$ & $\OO(uv)$ &&& \\ \hline
  $\OO(25s^2)$ & $\OO(uv+v^2)$ &&& \\ \hline
   $\OO(50s^2)$ & $\OO(uv+2v^2)$ &&& \\ \hline
   & $\{0\}$ & $\OO(y_1y_2)$&& \\ \hline
    & $\OO(u^2)$ & $\OO(y_1y_2+y_2^2)$&& \\ \hline
     & $\OO(2u^2)$ & $\OO(y_1y_2+2y_2^2)$&& \\ \hline
     && $\OO(\beta(x_1x_2x_3))$& $\{0\}$& \\ \hline
  && $\OO(\beta(x_1x_2x_3) + y_1^2)$& $\OO(z_1^2)$& \\ \hline
   && $\OO(\beta(x_1x_2x_3) + 2y_1^2)$& $\OO(2z_1^2)$& \\ \hline
     && $\OO(\beta(x_1x_2x_3) + 2y_1^2 +y_2^2)$& $\OO(2z_1^2+z_2^2)$& \\ \hline
    && $\OO(y_2y_3+ \beta(x_1x_2x_3))$& $\OO(z_1z_2)$&$\{0\}$ \\ \hline
   && $\OO(y_2y_3+ \beta(x_1x_2x_3) + y_1^2)$& &$\{\gamma^2 \}$ \\ \hline
    && $\OO(y_2y_3+ \beta(x_1x_2x_3) + 2y_1^2)$& &$\{2\gamma^2 \}$ \\ \hline
        && $\OO(y_2y_3+ \beta(x_1x_2x_3) + 3y_1^2)$& &$\{3\gamma^2 \}$ \\ \hline
            && $\OO(y_2y_3+ \beta(x_1x_2x_3) + 4y_1^2)$& &$\{4\gamma^2 \}$ \\ \hline
\end{tabular}}
\end{center}

    \bibliographystyle{alpha}
\bibliography{Categorical}

\def\cprime{$'$} \def\cprime{$'$}
\begin{thebibliography}{MnU18}

\bibitem[Lew68]{Lewis}
Gene Lewis.
\newblock The integral cohomology rings of groups of order {$p^{3}$}.
\newblock {\em Trans. Amer. Math. Soc.}, 132:501--529, 1968.

\bibitem[May18]{Maya}
Kevin Maya.
\newblock Clasificaci\'on de categor\'ias de fusi\'on punteadas de dimensi\'on
  $p^3$ hasta equivalencia morita.
\newblock {\em Master Thesis, Universidad del Norte}, 2018.

\bibitem[MnU18]{MunozUribe}
\'Alvaro Mu\~noz and Bernardo Uribe.
\newblock Classification of {P}ointed {F}usion {C}ategories of dimension 8 up
  to weak {M}orita equivalence.
\newblock {\em Comm. Algebra}, 46(9):3873--3888, 2018.

\bibitem[MS17]{MignardSchauenburg}
Micha\"el Mignard and Peter Schauenburg.
\newblock Morita equivalence of pointed fusion categories of small rank.
\newblock {\em arXiv:1708.06538}, 2017.

\bibitem[M{\"u}g03]{MugerI}
Michael M{\"u}ger.
\newblock From subfactors to categories and topology. {I}. {F}robenius algebras
  in and {M}orita equivalence of tensor categories.
\newblock {\em J. Pure Appl. Algebra}, 180(1-2):81--157, 2003.

\bibitem[Nai07]{Naidu}
Deepak Naidu.
\newblock Categorical {M}orita equivalence for group-theoretical categories.
\newblock {\em Comm. Algebra}, 35(11):3544--3565, 2007.

\bibitem[New72]{Newman}
Morris Newman.
\newblock {\em Integral matrices}.
\newblock Academic Press, New York-London, 1972.
\newblock Pure and Applied Mathematics, Vol. 45.

\bibitem[Uri17]{Uribe}
Bernardo Uribe.
\newblock On the classification of pointed fusion categories up to weak
  {M}orita equivalence.
\newblock {\em Pacific J. Math.}, 290(2):437--466, 2017.

\end{thebibliography}
   \end{document}